\documentclass{elsart}    

\usepackage{amsmath,amssymb,amsfonts,bm}
\usepackage{graphicx}
\usepackage{algorithm}
\usepackage{algorithmicx,algpseudocode}
\usepackage{cite}
\usepackage{color}
\let\theoremstyle\relax
\usepackage{amsthm}
\theoremstyle{plain}
\newtheorem{theorem}{Theorem}[section]
\newtheorem{lemma}[theorem]{Lemma}

\theoremstyle{definition}

\theoremstyle{remark}

\usepackage{wrapfig}

\begin{document}

\begin{frontmatter}

\title{Quantitative Error Analyses of Spectral Density Estimators Using Covariance Lags} 

\author[Guangyu]{Guangyu Wu}\ead{chinarustin@sjtu.edu.cn},    
\author[Anders]{Anders Lindquist}\ead{alq@kth.se}              

\address[Guangyu]{Department of Automation, Shanghai Jiao Tong University, Shanghai, China}  
\address[Anders]{Department of Automation and School of Mathematical Sciences, Shanghai Jiao Tong University, Shanghai, China}


\begin{abstract}                          
Spectral density estimation is a core problem of system identification, which is an important research area of system control and signal processing. There have been numerous results on the design of spectral density estimators. However to our best knowledge, quantitative error analyses of the spectral density estimation have not been proposed yet. In real practice, there are two main factors which induce errors in the spectral density estimation, including the external additive noise and the limited number of samples. In this paper, which is a very preliminary version, we first consider a univariate spectral density estimator using covariance lags. The estimation task is performed by a convex optimization scheme, and the covariance lags of the estimated spectral density are exactly as desired, which makes it possible for quantitative error analyses such as to derive tight error upper bounds. We analyze the errors induced by the two factors and propose upper and lower bounds for the errors. Then the results of the univariate spectral estimator are generalized to the multivariate one.
\end{abstract}

\end{frontmatter}

\section{Introduction}
\label{Introduction}
In this paper, we propose to analyze the errors of spectral density estimation by covariance lags. Spectral density estimation is a fundamental topic of system control and signal processing, where there have been numerous results, e.g. \cite{ andersson1998manual, gillberg2009frequency, byrnes1995complete, zorzi2015interpretation, you2022generalized, soderstrom1981identification, glover1974parametrizations, yuen2002spectral, pawitan1994nonparametric, 1996SPECTRAL}. Among these treatments, a widely known treatment is to use the covariance lags to estimate the spectral density by a convex optimization scheme. Proposed and advocated by Chris Byrnes, Tryphon Giorgiou, Anders Lindquist, and their collaborators and epigones \cite{byrnes2001finite, byrnes2003convex, pavon2006georgiou}, this type of spectral density estimator has drawn wide interests from both the academia and the industry, and is still a research focus of recent papers. The algorithm is a convex optimization, and the parameters of the model are proved to be diffeomorphic to the covariance lags, which ensures the existence and uniqueness of the optimal solution. The most significant advantage of this type of algorithm is that the covariance lags of the spectral density estimate are exactly as desired. 

However, the previous results are mostly focused on the design of the estimators. Even error analysis is of great significance for the designed estimators to be used in real scenarios, there have been few results on it. In a recent pioneering result by Bin Zhu and Mattia Zorzi \cite{zhu2023statistical}, a consistency analysis of the spectral density estimator was considered, which serves as a solid foundation of the error analysis for the spectral density estimator. However, to prove the consistency of the spectral estimator is not enough for analyzing the error due to the following two factors. First, we are not provided with infinite number of data samples for the spectral estimation, which means that the statistics of the estimator, including the covariance lags, are biased. Moreover, in quite a lot applications, we are provided with only a little amount of data samples for spectral density estimation. It makes the error of estimation not ignorable. Second, there always exists an external additive noise sequence to the original signal sequence to be estimated. The additive noise can even be malicious attacks \cite{pasqualetti2012cyber}, which may cause the spectral estimate to be severely biased from the true one. Due to these two reasons, we would like to investigate the errors of the statistics of the covaraince lags and their effects on the error of spectral estimation.

In this paper, we consider the error analyses of the spectral density estimator caused by the two factors mentioned above. Rather than analyzing the errors empirically by the simulation results, we propose theoretical error upper and lower bounds for the spectral estimator considering the existence of the two factors respectively. In doing this, there exist two main problems. The first one is to come up with a proper measure to describe the difference of the spectral estimate from the true one. The second problem is how to map the errors of the statistics of the covariance lags to those of the spectral estimate. 

The paper is organized as follows. In Section 2, we give a brief review of the univariate spectral density estimator. The main results are reviewed, however we will not go over all the detailed proofs throughout the derivation of the estimator. In the following Section 3, we propose an error upper bound, in the sense of the total variation distance, for the univariate spectral density estimator using the covariance lags with an additive noise sequence. We then propose both error upper and lower bounds of the spectral density estimator with limited number of samples in Section 4 and Section 5. In Section 6, we briefly review a generalization of the univariate spectral estimator to the multivariate case, and put forward error bounds for the multivariate spectral density estimator. A concluding remark is given in Section 7. 

\section{A brief review of the univariate spectral density estimator}
\label{Review}

In this section, we would like to first give a brief review of the spectral density estimator using the covariance lags by a convex optimization scheme proposed in \cite{byrnes2003convex}. The error analyses of the following two sections will be based on the estimator introduced in this section.

Denote $\mathcal{I} := [-\pi, \pi]$. A stationary stochastic process $y$ has a rational spectral density
$$
\Phi\left(e^{i \theta}\right)=\left|w\left(e^{i \theta}\right)\right|^2
$$
which is positive for all $\theta$. Its spectral density has a Fourier expansion
$$
\Phi\left(e^{i \theta}\right)=r_0+2 \sum_{k=1}^{\infty} r_k \cos k \theta,
$$
where the Fourier coefficients
\begin{equation}
r_k=\frac{1}{2 \pi} \int_{\mathcal{I}} e^{i k \theta} \Phi\left(e^{i \theta}\right) d \theta
\label{Covariancelag}
\end{equation}
are the covariance lags $r_k=\mathbb{E}[y(t+k) y(t)]$. The spectral density $\Phi(z)$ is analytic in an annulus containing the unit circle and has there the representation
$$
\Phi(z)=f(z)+f\left(z^{-1}\right),
$$
where $f$ is a rational function with all its poles and zeros in the open unit disc. Moreover, $\Phi\left(e^{i \theta}\right)=2 \operatorname{Re}\left\{f\left(e^{i \theta}\right)\right\}>0$ for all $\theta$, and therefore $f$ is a real function which maps $\{|z| \geq 0\}$ into the right half-plane $\{\operatorname{Re} z>0\}$; such a function is called positive real. For this to hold, the Toeplitz matrices
$$
T_n=\left[\begin{array}{cccc}
r_0 & r_1 & \cdots & r_n \\
r_1 & r_0 & \cdots & r_{n-1} \\
\vdots & \vdots & \ddots & \vdots \\
r_n & r_{n-1} & \cdots & r_0
\end{array}\right]
$$
must be positive definite for $n \in \mathbb{N}_0$, where $\mathbb{N}_0$ denotes the non-negative integers.

We first briefly paraphrase the results in \cite{byrnes2003convex}. We now refer the spectral density estimation problem treated in this paper to estimating $\Phi(z)$ by the covariance lags $r = \left(r_0, \cdots, r_n \right)$. 

Define the open convex cone $\mathfrak{P}_{+} \subset \mathbb{R}^{n}$ of sequences $q=$ $\left(q_1, q_2, \cdots, q_n\right)$ such that the corresponding generalized polynomial
\begin{equation}
Q(e^{i\theta})=\sum_{k=0}^n q_k e^{ik\theta}
\label{QxUni}
\end{equation}
is positive for all $ \theta \in \mathcal{I}$. For any choice of $P \in L_{+}^1([-\pi, \pi])$, the constrained optimization problem to minimize the functional
$$
\mathbb{I}_{P}(\Phi)=\int_{\mathcal{I}} P(e^{i\theta}) \log \frac{P(e^{i\theta})}{\Phi(e^{i\theta})} d \theta
$$
over $L_{+}^1([-\pi, \pi])$ subject to the constraints
$$
\int_{\mathcal{I}} e^{ik\theta} \Phi(e^{ik\theta}) d \theta=r_k, \quad k=0,1, \cdots, n,
$$
has unique solution, and it has the form
\begin{equation}
\Phi=\frac{P}{Q},
\label{FuncPhi}
\end{equation}
where $q \in \mathfrak{P}_{+}$ is the unique minimum of the strictly convex functional
\begin{equation}
\mathbb{J}_{P}(Q)=\langle r, q\rangle-\int_{\mathcal{I}} P \log Q d \theta.
\label{optimizationQ}
\end{equation}

In plain words, by the results in \cite{byrnes2003convex}, we have that a spectral density estimate in the form of \eqref{FuncPhi} can be uniquely determined by the covariance lags. Moreover, a very significant advantage of the spectral estimator proposed in \cite{byrnes2003convex} is that the covariance lags  are exactly as specified, which makes it possible for us to come up with error upper and lower bounds of it. In the following section, we will first consider analyzing the error of the spectral estimator introduced by an additive noise sequence.

\section{An error upper bound of univariate spectral density estimation with an additive noise sequence}

In this section, we will analyze the error of the spectral estimator induced by an additive noise sequence. Denote the stationary stochastic process corrupted with an additive noise sequence as $\tilde{y}=y+w$, where $w$ is also a stationary stochastic process and is independent of $y$.
Assume the first order moments $\mathbb{E}[y_{t}]=0$, and those of $w_{t}$ as $\mathbb{E}[w_{t}]=0$. Let the second order moments of $w_{t}$, namely its covariance lags, be $\mathbb{E}[w_{t+k} w_{t}]=r_{w, k}$. Then the covariance lags of $\tilde{y}_{t}$ read
$$
\begin{aligned}
& \mathbb{E}[\tilde{y}_{t+k} \tilde{y}_{t}] \\
= & \mathbb{E}[(y_{t+k}+w_{t+k})(y_{t}+w_{t})] \\
= & \mathbb{E}[y_{t+k} y_{t}]+\mathbb{E}[w_{t+k} w_{t}] \\
+ & \mathbb{E}[y_{t+k}] \mathbb{E}[w_{t}]+\mathbb{E}[y_{t}] \mathbb{E}[w_{t+k}] \\
= & r_k+r_{w, k}
\end{aligned}
$$

We define the spectral density of the stationary stochastic process $\tilde{y}$ as $\tilde{\Phi}(z)$. The problem now comes to selecting a proper metric to measure the difference between $\Phi(z)$ and $\tilde{\Phi}(z)$, and then propose an error upper bound in the sense of the metric.

We would first like to introduce some concepts from information theory. The total variation distance between the spectral density $\tilde{\Phi}$ with an additive noise sequence and the true density $\Phi$ is defined as follows:
\begin{equation}
    V(\tilde{\Phi}, \Phi) = \sup_{\theta} \left|\int_{\left[-\pi, \theta \right]}\left(\tilde{\Phi}\left(  e^{i\theta}\right) - \Phi\left( e^{i\theta} \right)\right) d \theta\right|
\label{lim}    
\end{equation}

In \cite{Aldo2003A}, Shannon-entropy is used to calculate the upper bound of the total variation distance. The Shannon-entropy \cite{shannon1948mathematical} modified for the spectral densities is defined as
$$
    H[\Phi] = - \int_{\mathcal{I}}\Phi\left( e^{i\theta} \right) \log \Phi\left( e^{i\theta} \right)d\theta.
$$

Then we introduce the concept of the maximum entropy distribution, which can be obtained by solving the following optimization problem which satisfies the  moment constraints
\begin{equation}
\begin{array}{ll}
\min_{\Phi} & -H(\Phi) \\
\text { s.t. } & \Phi(e^{i\theta}) \geq 0 \\
& \int_{\mathcal{I}} e^{ik\theta}\Phi(e^{i\theta}) d \theta=r_k\\ & \text { for } k=0, \cdots, n.
\end{array}
\label{OptBrevePhi}
\end{equation}

In the following part of this section, we will derive the form of the maximum entropy distribution subject to the moment constraints. Then we will present a formal proof that the derived spectral density is the maximum entropy distribution subject to given constraints.

The Lagragian of optimization \eqref{OptBrevePhi} can be written as
$$
\begin{aligned}
L(\Phi, \lambda) = & -H(\Phi)+\sum_{k=0}^n \lambda_k \left(\int_{\mathcal{I}} e^{ik\theta}\Phi(e^{i\theta}) d \theta - r_{k}\right)\\
+ & \lambda_{n+1}\Phi(e^{i\theta}).
\end{aligned}
$$

For the rest of the derivation, we will use the crude argument that we can think of the spectral density function $\Phi$ as an infinite-dimensional continuous vector with $\Phi\left( e^{i\theta} \right)$ as the value at each coordinate. Under this simplification, $\int_{\mathcal{I}} e^{ik\theta}\Phi(e^{i\theta}) d \theta$ is similar to $\sum_{\theta} e^{ik\theta}\Phi(e^{i\theta})$. We can then take derivative of $L(\Phi, \lambda)$ with respect to $\Phi(e^{i\theta}) \equiv \Phi_{\theta}$ treating $\Phi$ is a vector and all integrals as just summations.
$$
\frac{\partial L(\Phi, \lambda)}{\partial \Phi(e^{i\theta})}=\frac{\Phi(e^{i\theta})}{\Phi(e^{i\theta})}+\log \Phi(e^{i\theta})+\sum_{k=0}^n \lambda_k e^{ik\theta}+\lambda_{n+1}
$$
Setting $\frac{\partial L(\Phi, \lambda)}{\partial \Phi(e^{i\theta})}=0$ for all $\theta$, we have that the entropy maximizing distribution $\breve{\Phi}$ has the following form of function
$$
\breve{\Phi}(e^{i\theta})=\exp \left(-1-\sum_{k=0}^n \lambda_k e^{ik\theta}\right).
$$

We also note that the Shannon-entropy maximizing distribution satisfies the following equation
\begin{equation}
\begin{aligned}
H\left(\breve{\Phi}\right) = & -\int_{\mathcal{I}} \breve{\Phi}(e^{i\theta})\left(-1-\sum_{k=0}^n \lambda_k e^{ik\theta}\right)d\theta\\
= & 2\pi r_{0}+2\pi\sum_{k=0}^n \lambda_{k} r_k.
\end{aligned}
\label{HbrevePhi}
\end{equation}

Next we formally prove that $\breve{\Phi}$, as derived above, is indeed the maximum entropy distribution. We denote the Kullback-Leibler distance between the spectral densities $\Phi$ and $\breve{\Phi}$ as $KL\left( \Phi \| \breve{\Phi} \right)$, which is calculated by
$$
\begin{aligned}
    KL \left(\Phi \| \breve{\Phi}\right) & =  \int_{\mathcal{I}}\Phi(e^{i\theta}) \log \frac{\Phi(e^{i\theta})}{\breve{\Phi}(e^{i\theta})} d\theta \\
    & =- H\left [ \Phi \right ] + \sum_{k = 0}^{n} \lambda_{k}r_{k},
\end{aligned}
$$

\begin{lemma}
For all distributions $\Phi$ that satisfy the moment constraints, we have
$$
H\left(\breve{\Phi}\right) \geq H\left(\Phi\right)
$$
\end{lemma}
\begin{proof}
The Shannon-entropy of $\Phi$ reads
$$
\begin{aligned}
& H(\Phi) \\
= & -\int_{\mathcal{I}} \Phi(e^{i\theta}) \log \Phi(e^{i\theta}) \frac{\breve{\Phi}(e^{i\theta})}{\breve{\Phi}(e^{i\theta})}d\theta \\
= & -KL\left(\Phi \| \breve{\Phi}\right)-\int_{\mathcal{I}} \Phi(e^{i\theta}) \log \breve{\Phi}(e^{i\theta})d\theta \\
\leq & -\int_{\mathcal{I}} \Phi(e^{i\theta}) \log \breve{\Phi}(e^{i\theta}) \\
= & -\int_{\mathcal{I}} \Phi(e^{i\theta})\left(-1-\sum_{k=0}^n \lambda_k e^{ik\theta}\right)d\theta \\
= & 2\pi r_{0}+2\pi\sum_{k=0}^n \lambda_{k} r_k
\end{aligned}
$$

By \eqref{HbrevePhi}, we have
$$
\begin{aligned}
& H(\Phi) \\
= & -\int_{\mathcal{I}} \breve{\Phi}(e^{i\theta})\left(-1-\sum_{k=0}^n \lambda_k e^{ik\theta}\right)d\theta\\
= & -\int_{\mathcal{I}} \breve{\Phi}(e^{i\theta}) \log \breve{\Phi}(e^{i\theta})=H\left(\breve{\Phi}\right).
\end{aligned}
$$
which completes the proof to the lemma.
\end{proof}

With all the concepts introduced and lemma proved above, we now settle down to derive the error upper bound of the spectral density $\tilde{\Phi}$ of a stationary stochastic process $y$ with an additive noise sequence $w$.

We note that the Shannon-entropy in \eqref{lim} is upper bounded by
\begin{equation}
\begin{aligned}
& V(\tilde{\Phi}, \Phi)\\
\leq &\sup_{\theta} \left|\int_{\mathcal{I}}^{\theta}\left(\tilde{\Phi}\left(  e^{i\theta}\right) - \breve{\tilde{\Phi}}\left( e^{i\theta} \right)\right) d \theta\right|\\
+ & \sup_{\theta} \left|\int_{\mathcal{I}}^{\theta}\left(\breve{\tilde{\Phi}}\left(  e^{i\theta}\right) - \breve{\Phi}\left( e^{i\theta} \right)\right) d \theta\right|\\
+ & \sup_{\theta} \left|\int_{\mathcal{I}}^{\theta}\left(\breve{\Phi}\left(  e^{i\theta}\right) - \Phi\left( e^{i\theta} \right)\right) d \theta\right|\\
= & V(\tilde{\Phi}, \breve{\tilde{\Phi}}) + V(\breve{\tilde{\Phi}}, \breve{\Phi}) + V(\breve{\Phi}, \Phi).
\end{aligned}
\label{VTildePhiBound}
\end{equation}
where the Shannon-entropy maximizing distribution $\breve{\Phi}$ can be obtained by \eqref{OptBrevePhi}, and $\breve{\tilde{\Phi}}$ by
\begin{equation}
\begin{array}{ll}
\min_{\tilde{\Phi}} & -H(\tilde{\Phi}) \\
\text { s.t. } & \tilde{\Phi}(e^{i\theta}) \geq 0 \\
& \int_{\mathcal{I}} e^{ik\theta}\tilde{\Phi}(e^{i\theta}) d \theta=r_k+r_{w, k}\\
& \text { for } k=1, \cdots, n.
\end{array}
\label{OptBreveTildePhi}
\end{equation}

Define 
$$
    \Phi_{\text{m}} := \arg\max_{\rho \in \{ \Phi, \tilde{\Phi} \}} H\left( \rho \right).
$$

We note that the equality in \eqref{VTildePhiBound} is achieved if and only if $\Phi_{\text{m}} = \breve{\Phi}_{\text{m}}$.

Since we are able to obtain the analytic forms of function of both $\breve{\Phi}$ and $\breve{\tilde{\Phi}}$, we shall obtain $V(\breve{\tilde{\Phi}}, \breve{\Phi})$ by straightforward calculation. Now it remains to obtain $V(\tilde{\Phi}, \breve{\tilde{\Phi}})$ and $V(\breve{\Phi}, \Phi)$. By \cite{1970Correction, Aldo2003A}, we have
$$
\begin{aligned}
    V \left ( \tilde{\Phi}, \breve{\tilde{\Phi}}\right ) & \leq 3\left[-1+\left\{1+\frac{4}{9} KL \left(\tilde{\Phi} \| \breve{\tilde{\Phi}} \right)\right\}^{1 / 2}\right]^{1 / 2} \\
    & = 3\left[-1+\left\{1+\frac{4}{9} \left ( H\left [ \breve{\tilde{\Phi}} \right ] - H\left [ \tilde{\Phi} \right ] \right )\right\}^{1 / 2}\right]^{1 / 2}
\label{Vbound}
\end{aligned}
$$
and
$$
    V \left ( \breve{\Phi}, \Phi \right ) \leq 3\left[-1+\left\{1+\frac{4}{9} \left ( H\left [ \breve{\Phi} \right ] - H\left [ \Phi \right ] \right )\right\}^{1 / 2}\right]^{1 / 2}
$$

Then we obtain the upper bound of the error
\begin{equation}
\begin{aligned}
& V \left ( \tilde{\Phi}, \Phi \right ) \\
\leq & 3\left[-1+\left\{1+\frac{4}{9} \left ( H\left [ \breve{\tilde{\Phi}} \right ] - H\left [ \tilde{\Phi} \right ] \right )\right\}^{1 / 2}\right]^{1 / 2} \\
+ & V(\breve{\tilde{\Phi}}, \breve{\Phi})\\
+ & 3\left[-1+\left\{1+\frac{4}{9} \left ( H\left [ \breve{\Phi} \right ] - H\left [ \Phi \right ] \right )\right\}^{1 / 2}\right]^{1 / 2}.
\end{aligned}
\label{UpperBoundNoise}
\end{equation}

Moreover, by the reviewed spectral density estimator, $\tilde{\Phi}(e^{i\theta})$ is uniquely determined given $n$ and $\tilde{r}_{k}, k = 0, \cdots, n$ without bias. Therefore, we use $\tilde{\Phi}(e^{i\theta})$ to calculate the upper bound. We note that the spectral density $\tilde{\Phi}\left( e^{i\theta} \right)$ which has the form
\begin{equation}
\tilde{\Phi}\left( e^{i\theta} \right) = \frac{P(e^{i\theta})}{\tilde{Q}(e^{i\theta})}
\label{TildePhi}
\end{equation}
can be obtained by the optimization
\begin{equation}
\max\tilde{\mathbb{J}}_{P}(\tilde{Q})
\label{optimizationtildePQ}
\end{equation}
where the functional
$$
\begin{aligned}
\tilde{\mathbb{J}}_{P}( \tilde{Q}) =
& -r_0 \tilde{q}_0-r_1 \tilde{q}_1-\cdots-r_n \tilde{q}_n \\
& -\frac{1}{2 \pi} \int_{\mathcal{I}} P\left(e^{i \theta}\right) \log \frac{P\left(e^{i \theta}\right)}{\tilde{Q}\left(e^{i \theta}\right)} d \theta.
\end{aligned}
$$
Then $H\left [ \tilde{\Phi} \right ]$ can be calculated by the obtained $\tilde{\Phi}(e^{i\theta})$ in \eqref{TildePhi}. However, it is not feasible for us to obtain $\Phi(e^{i\theta})$ and the corresponding $H\left[ \Phi \right]$ in \eqref{UpperBoundNoise} since we are only provided with the sequences $y$ rather than $\Phi(e^{i\theta})$. A common mistake in doing it is to choose $\Phi\left( e^{i\theta} \right)$ in \eqref{VTildePhiBound} by optimization \eqref{optimizationQ}, which has the form of \eqref{FuncPhi}. However, it is not correct, since $\Phi$ obtained by \eqref{optimizationQ} is not the true spectral density. To be more specific, its first $n$ orders moments are identical to those of the true one, however its moments of orders $n+1$ to $+\infty$ are not necessarily equal to those of the true one. Instead, we adopt the bound of the Shannon-entropy in information theory \cite{polyanskiy2014lecture}. In our case, the upper bound can be interpreted as $H(\Phi) \geq 0$ with equality iff $\Phi$ is supported on finitely many discrete points within $\left[ -\pi, \pi \right]$. From an engineering perspective, we have that the non-negative $H(\Phi) = 0$ iff $\Phi$ has finitely many frequency components. We then have the following upper bound of error where there exists an additive noise sequence

\begin{equation}
\begin{aligned}
& V(\tilde{\Phi}, \Phi)\\
\leq & 3\left[-1+\left\{1+\frac{4}{9} \left ( H\left [ \breve{\tilde{\Phi}} \right ] - H\left [ \tilde{\Phi} \right ] \right )\right\}^{1 / 2}\right]^{1 / 2}\\
+ & V(\breve{\tilde{\Phi}}, \breve{\Phi})\\
+ & 3\left[-1+\left\{1+\frac{4}{9} \left ( H\left [ \breve{\Phi} \right ] \right )\right\}^{1 / 2}\right]^{1 / 2}.
\end{aligned}
\label{UpperboundNoise}
\end{equation}

In conclusion, we have proposed an error upper bound \eqref{UpperboundNoise}, in the sense of the total variation distance, for the spectral density estimate where there is an additive noise sequence to the original stochastic process. In the next section, we will turn to the quantitative error analysis of the density estimation with limited number of data samples.

\section{An error upper bound of univariate spectral density estimation with limited number of samples}

In real applications, we are not able to obtain infinite numbers of data samples. And in quite a lot of scenarios, we suffer from the lack of data samples. Assume the number of data samples we obtain to be $N$, the problem now comes to estimating finite windows of covariance lags from the sequence of observation
$$
y_0, y_1, y_2, \cdots, y_N
$$
of the process $\{y_{k} \mid k \in \mathbb{N}_{0}\}$. By assuming the stochastic process $y$ to be ergodic, the covariance lags can be estimated by
\begin{equation}
\tilde{r}_k=\frac{1}{N+1-k} \sum_{t=0}^{N-k} y_{t}y_{t+k}.
\end{equation}
However, due to the lack of data samples, it is possible for the estimates of covariance lags to be severely biased from the true ones. It is also mentioned in \cite{byrnes2002identifiability, byrnes2002identifiability1} that we can only estimate
$$
r_0, r_1, \cdots, r_n
$$
where $n \ll N$, with some precision. Hence the error analysis of spectral density estimation with limited number of data samples is of great significance both theoretically and empirically. 

We first investigate the probability of $\tilde{r}_{k}$ to fall within the interval $\left[ a_{k}, b_{k} \right]$, namely $\mathbb{P}\{ a_{k} \leq \tilde{r}_{k} \leq b_{k} \}$, considering three types of knowledge of $y_{0}, \cdots, y_{N}$.

For the first type of knowledge, we assume that the joint distribution of $y_{0}, \cdots, y_{N}$ is known prior. Then the probability can be obtained through direct calculation. However, even with a moderate amount of data samples, the calculation will be complicated.

For the second type of knowledge, we assume that only the marginal distributions of each $y_{t}y_{t+k}$ for $t = 0, \cdots, N - k$ are known. By the assumption that the stochastic process is stationary, all the marginal distributions are identical, i.e., 
$$
p\left( y_{0}y_{k} \right) = \cdots = p\left( y_{N-k}y_{N} \right).
$$

We note that
\begin{equation}
\begin{aligned}
& \mathbb{P}\{ a_{k} \leq \tilde{r}_{k} \leq b_{k} \}\\
= & \mathbb{P}\{ \left(N+1-k \right)a_{k} \leq \sum_{t=0}^{N-k} y_{t}y_{t+k}\\
& \leq \left(N+1-k \right)b_{k} \}\\
= & 1 - \mathbb{P}\{ \sum_{t=0}^{N-k} y_{t}y_{t+k} \leq \left(N+1-k \right)a_{k} \}\\
- & \mathbb{P}\{ \sum_{t=0}^{N-k} y_{t}y_{t+k} \geq \left(N+1-k \right)b_{k} \}\\
\leq & 1 - \mathbb{P}\{ y_{t}y_{t+k} \leq a_{k}, \forall 0 \leq t \leq N-k\}\\
- & \mathbb{P}\{ y_{t}y_{t+k} \geq b_{k}, \forall 0 \leq t \leq N-k\}\\
= & 1 - \left(\mathbb{P}\{ y_{t}y_{t+k} \leq a_{k}\}\right)^{N+1-k}\\
- & \left(\mathbb{P}\{ y_{t}y_{t+k} \geq b_{k}\}\right)^{N+1-k}
\end{aligned}
\label{arkb}
\end{equation}
Since the marginal distribution $p\left( y_{t}y_{t+k} \right)$ is known, we have
$$
\mathbb{P}\{ y_{t}y_{t+k} \leq a_{k}\} = \int_{-\infty}^{a_{k}} y_{t}y_{t+k} p\left( y_{t}y_{t+k} \right) dy_{t}y_{t+k}
$$
and
$$
\mathbb{P}\{ y_{t}y_{t+k} \geq b_{k}\} = \int_{b_{k}}^{+\infty} y_{t}y_{t+k} p\left( y_{t}y_{t+k} \right) dy_{t}y_{t+k}
$$

For the third type of knowledge, we only assume to know finite orders of power moments of $y_{t}y_{t+k}$, i.e., 
\begin{equation}
\int_{-\infty}^{+\infty} y_{t}^{i}y_{t+k}^{i} p\left( y_{t}y_{t+k} \right) dy_{t}y_{t+k}, 0 \leq i < +\infty.
\label{ytytkMoments}
\end{equation}
Even the knowledge is quite limited in this type which makes the problem a truncated Hamburger moment problem, we note that there are a series of research results on the bounds of the moment problem. These results make it feasible for us to derive upper bounds $\max \mathbb{P}\{y_{t}y_{t+k} \leq b_{k}\}$ and $\max \mathbb{P}\{y_{t}y_{t+k} \geq a_{k}\}$ given the moment constraints \eqref{ytytkMoments}. We can then obtain the lower bounds $\min \mathbb{P}\{y_{t}y_{t+k} \geq b_{k}\}$, $\min \mathbb{P}\{y_{t}y_{t+k} \leq a_{k}\}$, and then an upper bound of $\mathbb{P}\{ a_{k} \leq \tilde{r}_{k} \leq b_{k} \}$ in \eqref{arkb}. For example, the achievable upper bounds given the moment constraints derived by optimization schemes are proposed in \cite{bertsimas2005optimal}. 

We have proposed the upper and lower bounds of the covariance lags $\tilde{r}_{k}$ for $k = 0,\cdots, n$. Now the problem comes to putting forward an error upper bound of the spectral density estimate $\tilde{\Phi}$. We first consider the error upper bound in the sense of the total variation distance as in the previous section. 

We note that the problem of this section is different from the previous one. Given a specific additive noise sequence, $\tilde{\Phi}$ is determined in the previous section. However in this section, with limited number of data samples, $\tilde\Phi$ is not determinant, which makes $H\left( \tilde{\Phi} \right)$ in 
\eqref{UpperboundNoise} difficult to treat. We use the fact that $H\left[ \tilde{\Phi} \right]$ and obtain the following upper bound
\begin{equation}
\begin{aligned}
& V(\tilde{\Phi}, \Phi)\\
\leq & 3\left[-1+\left\{1+\frac{4}{9} \left ( H\left [ \breve{\tilde{\Phi}} \right ] \right )\right\}^{1 / 2}\right]^{1 / 2}\\
+ & 3\left[-1+\left\{1+\frac{4}{9} \left ( H\left [ \breve{\tilde{\Phi}} \right ] - H\left [ \breve{\Phi} \right ] \right )\right\}^{1 / 2}\right]^{1 / 2}\\
+ & 3\left[-1+\left\{1+\frac{4}{9} \left ( H\left [ \breve{\Phi} \right ] \right )\right\}^{1 / 2}\right]^{1 / 2}.
\end{aligned}
\label{UpperboundTotalNumber}
\end{equation}
We note that to maximize $V(\tilde{\Phi}, \Phi)$ is equivalent to maximize $H\left [ \breve{\tilde{\Phi}} \right ]$, i.e.,
$$
\arg\max_{\breve{\tilde{\Phi}}}V\left( \tilde{\Phi}, \Phi \right) = \arg\max_{\breve{\tilde{\Phi}}}H\left[ \breve{\tilde{\Phi}} \right].
$$
Then the problem is now obtaining the Shannon-entropy maximizing distribution, which can be formulated as the following optimization problem
\begin{equation}
\begin{array}{ll}
\max_{\tilde{\lambda}} & 2\pi \tilde{r}_{0}+2\pi\sum_{k=0}^n \tilde{\lambda}_{k} \tilde{r}_k \\
\text { s.t. } & a_{k} \leq \int_{\mathcal{I}} e^{ik\theta}\exp \left(-1-\sum_{k=0}^n \tilde{\lambda}_k e^{ik\theta}\right) d \theta\\
& = \tilde{r}_k \leq b_{k}\\ & \text { for } k=0, \cdots, n.
\end{array}
\label{Optiakbk}
\end{equation}
Denote the solution to this optimization problem as $\left( \tilde{\lambda}^{*}_{0}, \tilde{\lambda}^{*}_{1} \cdots, \tilde{\lambda}^{*}_{n} \right)$. Then the Shannon-entropy maximizing distribution reads
$$
\breve{\tilde{\Phi}}\left( e^{i\theta} \right) = \exp \left(-1-\sum_{k=0}^n \tilde{\lambda}^{*}_k e^{ik\theta}\right)
$$
and we can then obtain the corresponding $H\left[ \breve{\tilde{\Phi}} \right]$.

Moreover, by \eqref{Covariancelag}, we note that the covariance lags are indeed the projection of the spectral density $\Phi$ onto the orthonormal basis $\{ 1, e^{i\theta}, \cdots, e^{in\theta} \}$. Since there is no constraint on $\tilde{r}_{j}$ given all $\tilde{r}_{k}$ with $k \neq j$, we have that all $\tilde{r}_{k}$ for $k = 0, \cdots, n$ are independent, i.e., 

\begin{equation}
\begin{aligned}
& \mathbb{P}\{ a_{k} \leq \tilde{r}_{k} \leq b_{k}, \forall 0 \leq k \leq n \}\\
= & \prod_{k = 0}^{n}\mathbb{P}\{a_{k} \leq \tilde{r}_{k} \leq b_{k} \}\\
= & \prod_{k = 0}^{n}\left( 1 - \left(\mathbb{P}\{ y_{t}y_{t+k} \leq a_{k}\}\right)^{N+1-k} \right.\\
& \left. - \left(\mathbb{P}\{ y_{t}y_{t+k} \geq b_{k}\}\right)^{N+1-k} \right).
\end{aligned}
\end{equation}

In conclusion, the error upper bound can be interpreted as follows. The total variation distance between $\Phi$ and $\tilde{\Phi}$ has an upper bound \eqref{UpperboundTotalNumber} where $\breve{\tilde{\Phi}}\left( e^{i\theta} \right)$ is obtained by the optimization \eqref{Optiakbk}, with probability no greater than $\mathbb{P}\{ a_{k} \leq \tilde{r}_{k} \leq b_{k}, \forall 0 \leq k \leq n \}$. The upper bound reveals the fact that with the increase of the number of data samples $N$, it is more likely for the error in the sense of the total variation distance to fall within a specified compact interval.

\section{An error lower bound of univariate spectral density estimation with limited number of samples}

In the previous section, we proposed an error upper bound for spectral density estimation with limited number of data samples. However, in numerous scenarios, we would also like to derive an error lower bound to give us a sense of the least error that we shall have in our estimate. In this section, we will propose such an error lower bound of $\tilde{\Phi}$ in the sense of the Kullback-Leibler (KL) distance. Unlike the previous results where only the covariance lags are used to derive the upper bounds, the lower bound in the sense of the KL distance is directly related to the cepstral coefficients. We assume that the true spectral density $\Phi$ is known prior. The KL distance between the true spectral density and the density estimate by $N$ data samples reads
\begin{equation}
\begin{aligned}
    & KL \left(\Phi \| \tilde{\Phi}\right)\\ = & \int_{\mathcal{I}}\Phi(e^{i\theta}) \log \frac{\Phi(e^{i\theta})}{\tilde{\Phi}(e^{i\theta})} d\theta \\
    = & - H\left [ \Phi \right ] - \int_{\mathcal{I}}\Phi(e^{i\theta}) \log \tilde{\Phi}(e^{i\theta})d\theta
\end{aligned}
\label{KLPhitilde}
\end{equation}

Since $\Phi$ is known prior, it remains to treat the second term of \eqref{KLPhitilde}. We denote 
$$
\Phi(e^{i\theta}) := \sum_{k = 0}^{+\infty} \mu_{k} e^{ik\theta}.
$$

Another well-known way of representing the distribution of the stationary stochastic process is via the so-called cepstrum \cite{byrnes2002identifiability, byrnes2002identifiability1}
$$
\log \Phi\left(e^{i \theta}\right)=c_0+2 \sum_{k=1}^{\infty} c_k \cos k \theta
$$
The Fourier coefficients
$$
c_k=\frac{1}{2 \pi} \int_{\mathcal{I}} e^{i k \theta} \log \Phi\left(e^{i \theta}\right) d \theta
$$
are known as the cepstral coefficients.

With a proper choice of $n$, we have the following approximation
$$
\begin{aligned}
& -\int_{\mathcal{I}}\Phi(e^{i\theta}) \log \tilde{\Phi}(e^{i\theta})\\
= & -\sum_{k=0}^{+\infty} \mu_{k} \int_{\mathcal{I}}e^{ik\theta}\log \tilde{\Phi}(e^{i\theta}) d\theta\\
= & -\sum_{k=0}^{+\infty} 2\pi \mu_{k} \tilde{c}_{k} \approx -\sum_{k=0}^{n} 2\pi \mu_{k} \tilde{c}_{k}
\end{aligned}
$$
where $\mu_{k} \in \mathbb{R}_{+}$. Then the problem comes to deriving the lower bound of each $\tilde{c}_{k}$. In \cite{byrnes2002identifiability, byrnes2002identifiability1}, the cepstral coefficients are not obtained by some statistics of the data samples. Instead, they are tuned artificially to better fit the spectral density estimate, in the form of \eqref{FuncPhi}, to the data samples. However, we have the following inequality
\begin{equation}
\log \Phi(e^{i\theta}) \leq \Phi(e^{i\theta}) - 1.
\label{LogIneq}
\end{equation}
Therefore, since $\mu_{k}$ are all nonnegative, we have
$$
\begin{aligned}
& -\sum_{k=0}^{n} \mu_{k} \int_{\mathcal{I}}e^{ik\theta}\log \tilde{\Phi}(e^{i\theta}) d\theta\\
\geq & -\sum_{k=0}^{n} \mu_{k} \int_{\mathcal{I}}e^{ik\theta} \left(\tilde{\Phi}(e^{i\theta}) - 1\right) d\theta\\
= & -\sum_{k=0}^{n} \mu_{k}\tilde{r}_{k} 
\end{aligned}
$$
where
$$
\int_{\mathcal{I}} e^{ik\theta}d\theta = \int_{[-\pi, 0]}e^{ik\theta}d\theta + \int_{[0, \pi]} e^{ik\theta}d\theta = 0.
$$

Therefore, by \eqref{arkb} we have
$$
\begin{aligned}
-\sum_{k=0}^{n} \mu_{k} \int_{\mathcal{I}}e^{ik\theta}\log \tilde{\Phi}(e^{i\theta}) d\theta \geq -\sum_{k=0}^{n} \mu_{k}b_{k}\\
\end{aligned}
$$
with probability no greater than $\mathbb{P}\{ a_{k} \leq \tilde{r}_{k} \leq b_{k}, \forall 0 \leq k \leq n \}$.

Then we have the lower bound of error
\begin{equation}
KL \left(\Phi \| \tilde{\Phi}\right) \geq -\sum_{k=0}^{n} \mu_{k}b_{k} - H\left[ \Phi \right]
\label{KLlowerBound}
\end{equation}
with probability no greater than $\mathbb{P}\{ a_{k} \leq \tilde{r}_{k} \leq b_{k}, \forall 0 \leq k \leq n \}$.

Here we note that $-\sum_{k=0}^{n} \mu_{k}b_{k}$ needs to be nonnegative and greater than $H\left[ \Phi \right]$, or the error lower bound in \eqref{KLlowerBound} shall be negative, which is trivial since the Kullback-Leibler distance is always nonnegative.

\section{A brief review of the multivariate spectral density estimator and the corresponding error bounds}

In the previous sections, we have considered a univariate spectral density estimator using covariance lag by a convex optimization scheme. However, in quite some modern applications, e.g. 
image and signal processing, the spectral density to estimate is multiple dimensional. In these scenarios, the stationary stochastic process $\mathbf{y} = \left( y_{1}, \cdots, y_{d} \right)$ has multiple dimensions. A multiple dimensional spectral density estimator by a convex optimization scheme, of which the covariance lags are exactly as specified, is proposed in \cite{karlsson2016multidimensional}.

We first briefly review the results in \cite{karlsson2016multidimensional}, which is a generalization of the univariate spectral density estimator \cite{byrnes2003convex} to the multivariate case. Since \cite{karlsson2016multidimensional} treats the general multivariate moment problem, we paraphrase the results for the trigonometric moment problem. Define $\boldsymbol{\theta} := \left( \theta_{1}, \cdots, \theta_{d} \right) \in \mathcal{I}^{d}$. Let $\left\{\alpha_0, \cdots, \alpha_n\right\}$ be a set of trigonometric polynomials defined on $\mathcal{I}^d$, where
$$
\alpha_{k}(e^{i\boldsymbol{\theta}}) = \prod_{i=1}^{d}e^{i\alpha_{k, i}\theta_{i}}.
$$
The parameters $\alpha_{k, i} \in \mathbb{N}_{0}$ for $i = 1, \cdots, d$ are set such that the functions $\alpha_0, \alpha_1, \cdots, \alpha_n$ are linearly independent. Define the open convex cone $\mathfrak{P}_{+} \subset \mathbb{R}^n$ of sequences $p=$ $\left(p_1, p_2, \cdots, p_n\right)$ such that the corresponding generalized polynomial
$$
P(e^{i\boldsymbol{\theta}})=\sum_{k=0}^n p_k \alpha_k(e^{i\boldsymbol{\theta}})
$$
is positive for all $\boldsymbol{\theta}=\left(\theta_1, \cdots, \theta_d\right) \in \mathcal{I}^{d}$. Moreover, we denote by $\overline{\mathfrak{P}}_{+}$ its closure and by $\partial \mathfrak{P}_{+}$ its boundary $\overline{\mathfrak{P}}_{+} \backslash \mathfrak{P}_{+}$. We note that $P \equiv 0$ if and only if $p=0$, since $\alpha_1, \alpha_2, \cdots, \alpha_n$ are linearly independent. Denote $$
\mathfrak{R}_{+}=\left\{r \in \mathbb{R}^n \mid\langle r, p\rangle>0 \text { for all } p \in \overline{\mathfrak{P}}_{+} \backslash\{0\}\right\},
$$
where $\langle r, p\rangle$ is the inner product
$$
\langle r, p\rangle=\sum_{k=0}^n r_k p_k .
$$
Then we have the following theorem in \cite{karlsson2016multidimensional}. Denote 
$$
Q(e^{i\boldsymbol{\theta}}) =\sum_{k=0}^n q_k \alpha_k(e^{i\boldsymbol{\theta}}).
$$
Suppose that $(r, p) \in \mathfrak{R}_{+} \times \mathfrak{P}_{+}$ and the cone $\mathfrak{P}_{+}$ is nonempty and has the property
$$
\int_{\mathcal{I}^{d}} \frac{1}{Q} d \boldsymbol{\theta}=\infty \quad \text { for all } q \in \partial \mathfrak{P}_{+}.
$$
Then the optimization problem to maximize
$$
\mathbb{I}_P(\Phi)=\int_{\mathcal{I}^{d}} P(e^{i\boldsymbol{\theta}}) \log \Phi(e^{i\boldsymbol{\theta}}) d \boldsymbol{\theta}
$$
over all $\Phi \in L_{1}^{+}\left(\mathcal{I}^{d}\right)$ satisfying the moment condition
$$
\int_{\mathcal{I}^{d}} \alpha_{k}(e^{i\boldsymbol{\theta}}) \Phi(e^{i\boldsymbol{\theta}}) d \boldsymbol{\theta}=r_{k}, \ \text{for} \ k = 0, \cdots, n
$$
has a unique solution
$$
\Phi=\frac{P}{Q},
$$
where $Q$ is the unique minimizer of 
$$
\mathbb{J}_{P}^{r}(Q)=\langle r, q\rangle-\int_{\mathcal{I}^{d}} P \log Q d \boldsymbol{\theta}.
$$

Provided with the result in \cite{karlsson2016multidimensional}, we now settle down to carry out quantitative analyses of the multivariate spectral density estimator. 

The Shannon-entropy maximizing distribution for the multivariate case, namely $\breve{\Phi}$, can be obtained by the following optimization

\begin{equation}
\begin{array}{ll}
\min_{\Phi} & -H(\Phi) \\
\text { s.t. } & \Phi(e^{i\boldsymbol{\theta}}) \geq 0 \\
& \int_{\mathcal{I}^{d}} \alpha_{k}(e^{i\boldsymbol{\theta}})\Phi(e^{i\boldsymbol{\theta}}) d \boldsymbol{\theta}=r_k\\ & \text { for } k=0, \cdots, n.
\end{array}
\label{OptBrevePhiMulti}
\end{equation}

Then we can form the Lagragian as
$$
\begin{aligned}
L(\Phi, \lambda) = & -H(\Phi) + \lambda_{n+1}\Phi(e^{i\boldsymbol{\theta}})\\
+ & \sum_{k=0}^n \lambda_k \left(\int_{\mathcal{I}^{d}}\alpha_{k} \left(e^{i\boldsymbol{\theta}}\right)\Phi(e^{i\boldsymbol{\theta}}) d \boldsymbol{\theta} - r_{k}\right).
\end{aligned}
$$

Following a similar treatment as that in Section 3, we obtain the multivariate Shannon-entropy maximizing distribution
$$
\breve{\Phi}(e^{i\boldsymbol{\theta}})=\exp \left(-1-\sum_{k=0}^n \lambda_k \alpha_{k}\left(e^{i\boldsymbol{\theta}}\right)\right)
$$
by taking $\frac{\partial L(\Phi, \lambda)}{\partial \Phi(e^{i\boldsymbol{\theta}})}=0$.

Denote the stationary stochastic process corrupted with an additive noise sequence as $\tilde{\mathbf{y}}=\mathbf{y}+\mathbf{w}$, where $\mathbf{w}$ is also a stationary stochastic process and is independent of $\mathbf{y}$. The $i_\text{th}$ dimension of $\tilde{\mathbf{y}}$ at time step $t$ is denoted as $\mathbf{y}_{t, i}$. By assuming $\mathbb{E}[y_{t}] = \mathbb{E}[w_{t}] = 0$, we shall write the covariance lags of the multivariate spectral density as 
$$
\mathbb{E}[\tilde{\mathbf{y}}_{t+k} \tilde{\mathbf{y}}_{t}] = \mathbb{E}[\mathbf{y}_{t+k} \mathbf{y}_{t}]+\mathbb{E}[\mathbf{w}_{t+k} \mathbf{w}_{t}] = r_k+r_{w, k} = \tilde{r}_k.
$$

Then the Shannon-entropy maximizing distribution of the noise-corrupted stochastic process, namely $\breve{\tilde{\Phi}}$, can be obtained by
\begin{equation}
\begin{array}{ll}
\min_{\tilde{\Phi}} & -H(\tilde{\Phi}) \\
\text { s.t. } & \tilde{\Phi}(e^{i\boldsymbol{\theta}}) \geq 0 \\
& \int_{\mathcal{I}^{d}} \alpha_{k}(e^{i\boldsymbol{\theta}})\tilde{\Phi}(e^{i\boldsymbol{\theta}}) d \boldsymbol{\theta}=\tilde{r}_k\\
& \text { for } k=1, \cdots, n.
\end{array}
\label{OptBreveTildePhiMulti}
\end{equation}

Moreover, we assume that the true spectral density $\Phi(e^{i\boldsymbol{\theta}})$ is not known prior, except for the covariance lags. The spectral density estimate of $\tilde{\mathbf{y}}$, which has the form 
\begin{equation}
\tilde{\Phi}\left( e^{i\boldsymbol{\theta}} \right) = \frac{P(e^{i\boldsymbol{\theta}})}{\tilde{Q}(e^{i\boldsymbol{\theta}})},
\label{TildePhiMulti}
\end{equation}
can be obtained by minimizing
$$
\mathbb{J}_{P}^{\tilde{r}}(\tilde{Q}) = \langle \tilde{r}, \tilde{q}\rangle-\int_{\mathcal{I}^{d}} P \log \tilde{Q} d \boldsymbol{\theta}.
$$

In conclusion, with $\tilde{\Phi}$ obtained in \eqref{TildePhiMulti} and the multivariate Shannon-entropy maximizing distributions $\breve{\Phi}$ and $\breve{\tilde{\Phi}}$ obtained by optimizations \eqref{OptBrevePhiMulti} and \eqref{OptBreveTildePhiMulti}, an error upper bound of the multivariate stochastic process $\tilde{\mathbf{y}}$ corrupted with an additive sequence in the sense of the total variation distance can be calculated by \eqref{UpperboundNoise}.

In the following part of this section, we will analyze the error of estimation for the multivariate spectral densities. Similar to the results in Section 4, the covariance lags of the multivariate spectral density $\Phi(e^{i\boldsymbol{\theta}})$ can be estimated by
\begin{equation}
\tilde{r}_k=\frac{1}{N+1-k} \sum_{t=0}^{N-k} \prod_{i = 1}^{d}\mathbf{y}_{t, i}\mathbf{y}_{t+\alpha_{k, i}, i}.
\end{equation} by assuming the multivariate stochastic process $\mathbf{y}$ to be ergodic. We derive the probability of $\tilde{r}_{k}$ to fall within the interval $\left[ a, b \right]$, namely $\mathbb{P}\{ a \leq \tilde{r}_{k} \leq b \}$, considering two types of knowledge of $\mathbf{y}_{0}, \cdots, \mathbf{y}_{N}$. If the joint distribution of $\mathbf{y}_{0}, \cdots, \mathbf{y}_{N}$ is known, we shall obtain $\mathbb{P}\{ a_{k} \leq \tilde{r}_{k} \leq b_{k} \}$ by direct calculation. For the second type of knowledge, we assume that only the marginal distributions of each $\mathbf{y}_{t}\mathbf{y}_{t+k}$ for $t = 0, \cdots, N - k$ are known. Similar to \eqref{arkb}, we have that
\begin{equation}
\begin{aligned}
& \mathbb{P}\{ a_{k} \leq \tilde{r}_{k} \leq b_{k} \}\\
= & \mathbb{P}\{ \left(N+1-k \right)a_{k} \leq \sum_{t=0}^{N-k} \prod_{i = 1}^{d}\mathbf{y}_{t, i}\mathbf{y}_{t+\alpha_{k, i}, i}\\
& \leq \left(N+1-k \right)b_{k} \}\\
= & 1 - \mathbb{P}\{ \sum_{t=0}^{N-k} \prod_{i = 1}^{d}\mathbf{y}_{t, i}\mathbf{y}_{t+\alpha_{k, i}, i} \leq \left(N+1-k \right)a_{k} \}\\
- & \mathbb{P}\{ \sum_{t=0}^{N-k} \prod_{i = 1}^{d}\mathbf{y}_{t, i}\mathbf{y}_{t+\alpha_{k, i}, i} \geq \left(N+1-k \right)b_{k} \}\\
\leq & 1 - \mathbb{P}\{ \mathbf{y}_{t, i}\mathbf{y}_{t+\alpha_{k, i}, i} \leq a_{k}, \forall 0 \leq t \leq N-k\}\\
- & \mathbb{P}\{ \mathbf{y}_{t, i}\mathbf{y}_{t+\alpha_{k, i}, i} \geq b_{k}, \forall 0 \leq t \leq N-k\}\\
= & 1 - \left(\mathbb{P}\{ \mathbf{y}_{t, i}\mathbf{y}_{t+\alpha_{k, i}, i} \leq a_{k}\}\right)^{N+1-k}\\
- & \left(\mathbb{P}\{ \mathbf{y}_{t, i}\mathbf{y}_{t+\alpha_{k, i}, i} \geq b_{k}\}\right)^{N+1-k}.
\end{aligned}
\label{arkbMulti}
\end{equation}
Following the treatment in Section 4, we first obtain the Shannon-entropy maximizing distribution $\breve{\tilde{\Phi}}$ by the following optimization
\begin{equation}
\begin{array}{ll}
\tilde{\lambda}^{*}_{k} = &\arg\max_{\tilde{\lambda}} \left( 2\pi \tilde{r}_{0}+2\pi\sum_{k=0}^n \tilde{\lambda}_{k} \tilde{r}_k \right) \\
\text { s.t. } & a_{k} \leq \int_{\mathcal{I}^{d}} \alpha_{k}(e^{i\boldsymbol{\theta}})\exp \left(-1-\sum \tilde{\lambda}_k \alpha_{k}(e^{i\boldsymbol{\theta}})\right) d \boldsymbol{\theta}\\
& = \tilde{r}_k \leq b_{k}\\ & \text { for } k=0, \cdots, n.
\end{array}
\label{OptiakbkMulti}
\end{equation}
Then we have
$$
\breve{\tilde{\Phi}}(e^{i\boldsymbol{\theta}})=\exp \left(-1-\sum_{k=0}^n \tilde{\lambda}^{*}_{k} \alpha_{k}\left(e^{i\boldsymbol{\theta}}\right)\right).
$$

Since $\breve{\Phi}$ can be obtained by \eqref{OptBrevePhiMulti}, we shall calculate $H(\breve{\Phi})$. With $H(\breve{\Phi})$ and $H(\breve{\tilde{\Phi}})$ both known, we shall obtain the error upper bound in the sense of the total variation distance for the multivariate spectral density estimator with $N$ data samples, namely $V(\tilde{\Phi}, \Phi)$, by \eqref{UpperboundTotalNumber}. Therefore, we can conclude that the total variation distance between the multivariate $\Phi$ and $\tilde{\Phi}$ has an upper bound \eqref{UpperboundTotalNumber} where $\breve{\tilde{\Phi}}\left( e^{i\theta} \right)$ is obtained by the optimization \eqref{OptiakbkMulti}, with probability no greater than $\mathbb{P}\{ a_{k} \leq \tilde{r}_{k} \leq b_{k}, \forall 0 \leq k \leq n \}$. Since $\alpha_{1}(e^{i\boldsymbol{\theta}}), \cdots, \alpha_{n}(e^{i\boldsymbol{\theta}})$ are linear independent, we have that $\tilde{r}_{1}, \cdots, \tilde{r}_{n}$ are independent. Therefore we have
\begin{equation}
\begin{aligned}
& \mathbb{P}\{ a_{k} \leq \tilde{r}_{k} \leq b_{k}, \forall 0 \leq k \leq n \}\\
= & \prod_{k = 0}^{n}\mathbb{P}\{a_{k} \leq \tilde{r}_{k} \leq b_{k} \}\\
= & \prod_{k = 0}^{n}\left( 1 - \left(\mathbb{P}\{ \mathbf{y}_{t, i}\mathbf{y}_{t+\alpha_{k, i}, i} \leq a_{k}\}\right)^{N+1-k} \right.\\
& \left. - \left(\mathbb{P}\{ \mathbf{y}_{t, i}\mathbf{y}_{t+\alpha_{k, i}, i} \geq b_{k}\}\right)^{N+1-k} \right).
\end{aligned}
\label{ProbUpperBound}
\end{equation}

At last, we would put forward a lower bound of error for the multivariate spectral density estimator. For the multivariate case, the error of estimation in the sense of the Kullback-Leibler distance reads
\begin{equation}
\begin{aligned}
    & KL \left(\Phi \| \tilde{\Phi}\right)\\
    = & \int_{\mathcal{I}^{d}}\Phi(e^{i\boldsymbol{\theta}}) \log \frac{\Phi(e^{i\boldsymbol{\theta}})}{\tilde{\Phi}(e^{i\boldsymbol{\theta}})} d\boldsymbol{\theta} \\
    = & - H\left [ \Phi \right ] - \int_{\mathcal{I}^{d}}\Phi(e^{i\boldsymbol{\theta}}) \log \tilde{\Phi}(e^{i\boldsymbol{\theta}})d\boldsymbol{\theta}.
\label{KLMulti}
\end{aligned}
\end{equation}

We note that with a proper choice of $n$ which is large enough, the projection of $\Phi(e^{i\boldsymbol{\theta}})$ onto the vector space spanned by the basis $\left\{ \alpha_{1}(e^{i\boldsymbol{\theta}}), \cdots, \alpha_{n}(e^{i\boldsymbol{\theta}}) \right\}$ shall be close to $\Phi(e^{i\boldsymbol{\theta}})$ itself, i.e.,
\begin{equation}
\Phi(e^{i\boldsymbol{\theta}}) \approx \sum_{k = 1}^{n} \mu_{k} \alpha_{k}(e^{i\boldsymbol{\theta}}) \quad \text{when} \ n \gg 1,
\label{PhiMultiApprox}
\end{equation}
where $\mu_{k} \geq 0$. 

Moreover, we note that the dimension of the vector space is $n$, since $\alpha_{1}(e^{i\boldsymbol{\theta}}), \cdots, \alpha_{n}(e^{i\boldsymbol{\theta}})$ are linear independent. By \eqref{PhiMultiApprox}, we can write \eqref{KLMulti} as
\begin{equation*}
\begin{aligned}
    & KL \left(\Phi \| \tilde{\Phi}\right)\\
    = & - H\left [ \Phi \right ] - \sum_{k = 1}^{n}\mu_{k}\int_{\mathcal{I}^{d}} \alpha_{k}(e^{i\boldsymbol{\theta}}) \log \tilde{\Phi}(e^{i\boldsymbol{\theta}})d\boldsymbol{\theta}\\
    = & - H\left [ \Phi \right ] -\sum_{k = 1}^{n}\mu_{k}\tilde{c}_{k} .
\end{aligned}
\end{equation*}

Still by the log inequality \eqref{LogIneq}, we have
$$
\begin{aligned}
& -\sum_{k=1}^{n} \mu_{k} \int_{\mathcal{I}^{d}} \alpha_{k}(e^{i\boldsymbol{\theta}}) \log \tilde{\Phi}(e^{i\boldsymbol{\theta}})d\boldsymbol{\theta}\\
\geq & -\sum_{k=0}^{n} \mu_{k} \int_{\mathcal{I}^{d}} \alpha_{k}(e^{i\boldsymbol{\theta}}) \left(\tilde{\Phi}(e^{i\boldsymbol{\theta}}) - 1\right)d\boldsymbol{\theta}\\
= & -\sum_{k=0}^{n} \mu_{k}\tilde{r}_{k}.
\end{aligned}
$$

Therefore, we have
$$
\begin{aligned}
-\sum_{k=1}^{n} \mu_{k} \int_{\mathcal{I}^{d}} \alpha_{k}(e^{i\boldsymbol{\theta}}) \log \tilde{\Phi}(e^{i\boldsymbol{\theta}})d\boldsymbol{\theta} \geq -\sum_{k=1}^{n} \mu_{k}b_{k}.
\end{aligned}
$$
with probability no greater than $\mathbb{P}\{ a_{k} \leq \tilde{r}_{k} \leq b_{k}, \forall 0 \leq k \leq n \}$ calculated by \eqref{ProbUpperBound}. 

In conclusion, we obtain the lower bound of error in the sense of the Kullback-Leibler distance for the multivariate spectral estimator, which reads
\begin{equation}
KL \left(\Phi \| \tilde{\Phi}\right) \geq -\sum_{k=1}^{n} \mu_{k}b_{k} - H\left[ \Phi \right]
\label{KLlowerBoundMulti}
\end{equation}
with probability no greater than $\mathbb{P}\{ a_{k} \leq \tilde{r}_{k} \leq b_{k}, \forall 0 \leq k \leq n \}$ in \eqref{ProbUpperBound}. Similar to the result in Section 5, we need to select $b_{k}$ properly to ensure that the r.h.s. of \eqref{KLlowerBoundMulti} is positive.

\section{A concluding remark}
Quantitative error analyses of the spectral density estimation is of great significance for better understanding of the estimation algorithm and for the application of the estimation algorithm to real scenarios. However, the problem is quite difficult since the conventional estimators are usually obtained by optimization to drive their statistics to as close as the desired ones. The statistics are not exactly the ones desired, which makes it not feasible for us to propose quantitative bounds for the estimators. Proposed by Chris Byrnes, Tryphon Giorgiou, Anders Lindquist, an estimator using the covariance lags by a convex optimization scheme is able to satisfy the desired statistics without bias. This great property makes it feasible for us to analyze the errors of spectral estimation quantitatively. In this paper, we consider two typical factors which introduce errors to estimation, namely additive noise and limited number of data samples, and analyze the errors they introduce. We propose an error upper bound for the univariate spectral density estimator with an additive noise sequence. And we propose both an upper and a lower bound for the estimator with limited number of data samples. The results of the univariate estimator are then generalized to the multivariate case.

\bibliographystyle{plain}
\bibliography{autosam}

\begin{thebibliography}{10}

\bibitem{andersson1998manual}
Lennart Andersson, Ulf J{\"o}nsson, Karl~Henrik Johansson, and Johan Bengtsson.
\newblock A manual for system identification.
\newblock {\em Laboratory Exercises in System Identification. KF Sigma i Lund
  AB. Department of Automatic Control, Lund Institute of Technology, Box}, 118,
  1998.

\bibitem{bertsimas2005optimal}
Dimitris Bertsimas and Ioana Popescu.
\newblock Optimal inequalities in probability theory: A convex optimization
  approach.
\newblock {\em SIAM Journal on Optimization}, 15(3):780--804, 2005.

\bibitem{byrnes2002identifiability1}
Christopher~I Byrnes, Per Enqvist, and Anders Lindquist.
\newblock Identifiability and well-posedness of shaping-filter
  parameterizations: A global analysis approach.
\newblock {\em SIAM journal on control and optimization}, 41(1):23--59, 2002.

\bibitem{byrnes2002identifiability}
Christopher~I Byrnes, Per Enqvist, and Anders Lindquist.
\newblock Identifiability of shaping filters from covariance lags, cepstral
  windows and markov parameters.
\newblock In {\em Proceedings of the 41st IEEE Conference on Decision and
  Control, 2002.}, volume~1, pages 246--251. IEEE, 2002.

\bibitem{byrnes2001finite}
Christopher~I Byrnes, Sergei~V Gusev, and Anders Lindquist.
\newblock From finite covariance windows to modeling filters: A convex
  optimization approach.
\newblock {\em SIAM review}, 43(4):645--675, 2001.

\bibitem{byrnes2003convex}
Christopher~I Byrnes and Anders Lindquist.
\newblock {\em A convex optimization approach to generalized moment problems}.
\newblock Springer, 2003.

\bibitem{byrnes1995complete}
Christopher~I Byrnes, Anders Lindquist, Sergei~V Gusev, and Alexey~S Matveev.
\newblock A complete parameterization of all positive rational extensions of a
  covariance sequence.
\newblock {\em IEEE Transactions on Automatic Control}, 40(11):1841--1857,
  1995.

\bibitem{gillberg2009frequency}
Jonas Gillberg and Lennart Ljung.
\newblock Frequency-domain identification of continuous-time arma models from
  sampled data.
\newblock {\em Automatica}, 45(6):1371--1378, 2009.

\bibitem{glover1974parametrizations}
Keith Glover and Jan Willems.
\newblock Parametrizations of linear dynamical systems: Canonical forms and
  identifiability.
\newblock {\em IEEE Transactions on Automatic Control}, 19(6):640--646, 1974.

\bibitem{karlsson2016multidimensional}
Johan Karlsson, Anders Lindquist, and Axel Ringh.
\newblock The multidimensional moment problem with complexity constraint.
\newblock {\em Integral equations and operator theory}, 84(3):395--418, 2016.

\bibitem{1970Correction}
S.~Kullback.
\newblock Correction to a lower bound for discrimination information in terms
  of variation.
\newblock {\em IEEE Transactions on Information Theory}, 16(5):652--652, 1970.

\bibitem{1996SPECTRAL}
M.~H. Neumann.
\newblock Spectral density estimation via nonlinear wavelet methods for
  stationary non-gaussian time series.
\newblock {\em Journal of Time Series Analysis}, 1996.

\bibitem{pasqualetti2012cyber}
Fabio Pasqualetti, Florian D{\"o}rfler, and Francesco Bullo.
\newblock Cyber-physical security via geometric control: Distributed monitoring
  and malicious attacks.
\newblock In {\em 2012 IEEE 51st IEEE Conference on Decision and Control
  (CDC)}, pages 3418--3425. IEEE, 2012.

\bibitem{pavon2006georgiou}
Michele Pavon and Augusto Ferrante.
\newblock On the georgiou-lindquist approach to constrained kullback-leibler
  approximation of spectral densities.
\newblock {\em IEEE transactions on Automatic Control}, 51(4):639--644, 2006.

\bibitem{pawitan1994nonparametric}
Yudi Pawitan and Finbarr O'sullivan.
\newblock Nonparametric spectral density estimation using penalized whittle
  likelihood.
\newblock {\em Journal of the American Statistical Association},
  89(426):600--610, 1994.

\bibitem{polyanskiy2014lecture}
Yury Polyanskiy and Yihong Wu.
\newblock Lecture notes on information theory.
\newblock {\em Lecture Notes for ECE563 (UIUC) and}, 6(2012-2016):7, 2014.

\bibitem{shannon1948mathematical}
Claude~Elwood Shannon.
\newblock A mathematical theory of communication.
\newblock {\em The Bell system technical journal}, 27(3):379--423, 1948.

\bibitem{soderstrom1981identification}
Torsten S{\"o}derstr{\"o}m.
\newblock Identification of stochastic linear systems in presence of input
  noise.
\newblock {\em Automatica}, 17(5):713--725, 1981.

\bibitem{Aldo2003A}
Aldo Tagliani.
\newblock A note on proximity of distributions in terms of coinciding moments.
\newblock {\em Applied Mathematics and Computation}, 145(2-3):195--203, 2003.

\bibitem{you2022generalized}
Junyao You, Chengpu Yu, Jian Sun, and Jie Chen.
\newblock Generalized maximum entropy based identification of graphical arma
  models.
\newblock {\em Automatica}, 141:110319, 2022.

\bibitem{yuen2002spectral}
Ka-Veng Yuen, Lambros~S Katafygiotis, and James~L Beck.
\newblock Spectral density estimation of stochastic vector processes.
\newblock {\em Probabilistic Engineering Mechanics}, 17(3):265--272, 2002.

\bibitem{zhu2023statistical}
Bin Zhu and Mattia Zorzi.
\newblock On the statistical consistency of a generalized cepstral estimator.
\newblock {\em arXiv preprint arXiv:2301.06784}, 2023.

\bibitem{zorzi2015interpretation}
Mattia Zorzi.
\newblock An interpretation of the dual problem of the three-like approaches.
\newblock {\em Automatica}, 62:87--92, 2015.

\end{thebibliography}

\begin{wrapfigure}{l}{20mm} 
\includegraphics[width=1in,height=1.25in,clip,keepaspectratio]{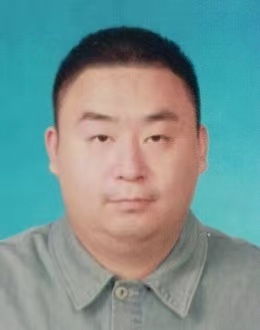}
  \end{wrapfigure}\par
  \textbf{Guangyu Wu} received the B.E. degree from Northwestern Polytechnical University, Xi’an, China, in 2013, and two M.S. degrees, one in control science and engineering from Shanghai Jiao Tong University, Shanghai, China, in 2016, and the other in electrical engineering from the University of Notre Dame, South Bend, USA, in 2018. 

He is currently pursuing the Ph.D. degree at Shanghai Jiao Tong University. His research interests are the moment problem and its applications to stochastic filtering, density steering, system identification and statistics. 

\begin{wrapfigure}{l}{25mm} 
\includegraphics[width=1in,height=1.25in,clip,keepaspectratio]{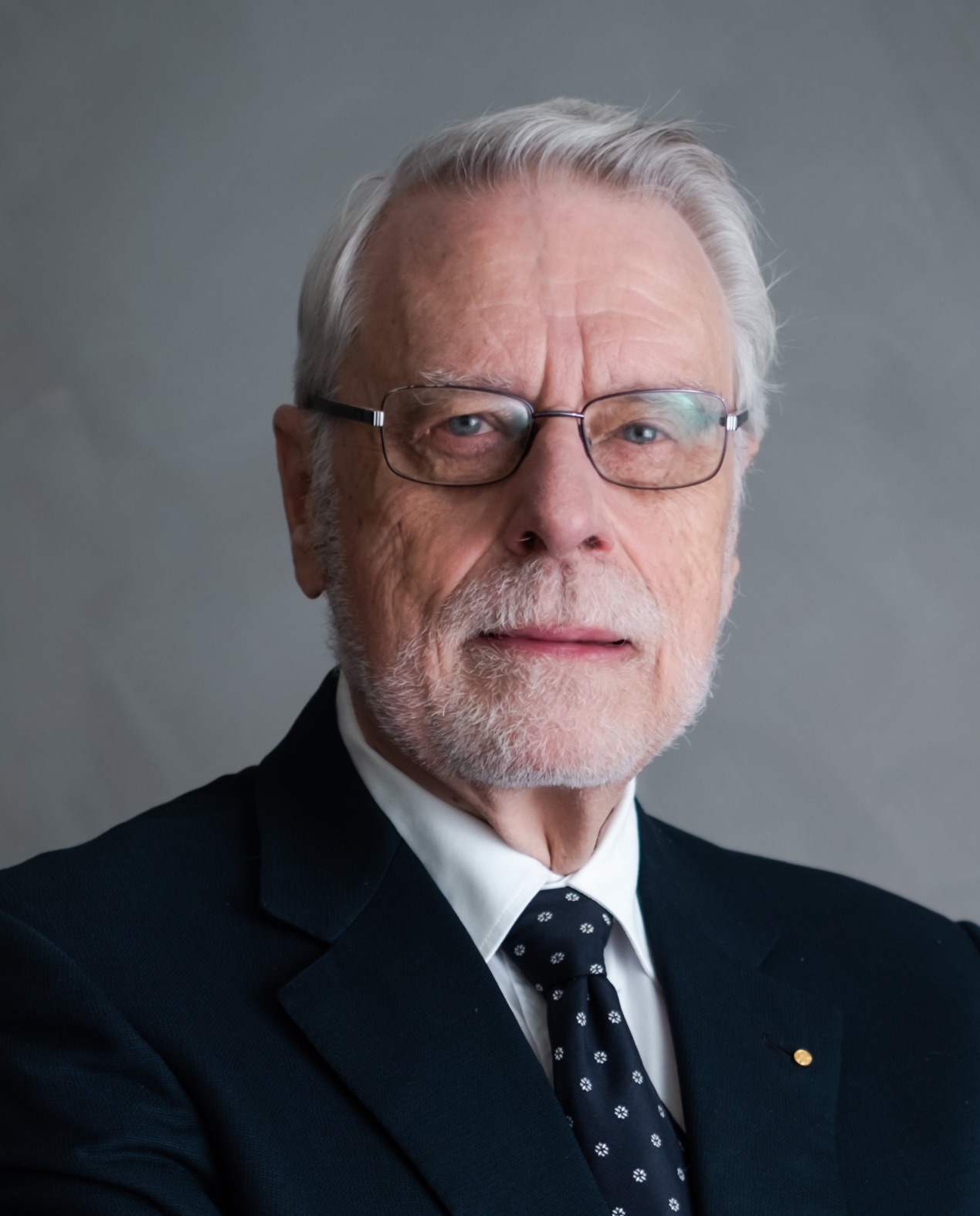}
  \end{wrapfigure}\par
  \textbf{Anders Lindquist} received the Ph.D. degree in optimization and systems theory from the Royal Institute of Technology (KTH), Stockholm, Sweden, in 1972, an honorary doctorate (Doctor Scientiarum Honoris Causa) from Technion (Israel Institute of Technology) in 2010 and Doctor Jubilaris from KTH in 2022.

He is currently a Zhiyuan Chair Professor at Shanghai Jiao Tong University, China, and Professor Emeritus at the Royal Institute of Technology (KTH), Stockholm, Sweden. Before that he had a full academic career in the United States, after which he was appointed to the Chair of Optimization and Systems at KTH.
Dr. Lindquist is a Member of the Royal Swedish Academy of Engineering Sciences, a Foreign Member of the Chinese Academy of Sciences, a Foreign Member of the Russian Academy of Natural Sciences, a Member of Academia Europaea (Academy of Europe), an Honorary Member the Hungarian Operations Research Society, a Fellow of SIAM, and a Fellow of IFAC. He received the 2003 George S. Axelby Outstanding Paper Award, the 2009 Reid Prize in Mathematics from SIAM, and the 2020 IEEE Control Systems Award, the
IEEE field award in Systems and Control.

\end{document}